\documentclass[12pt]{article}
\usepackage{amssymb, amsfonts, amsmath, amscd, theorem}
\usepackage[all]{xy}

\numberwithin{equation}{section}

\newcommand{\Aw}{\operatorname{A}^\lambda_{q,w}}
\newcommand{\Ae}{\operatorname{A}^\lambda_{q,e}}
\newcommand{\tD}{\operatorname{\widetilde{\mathcal{D}}_q}}
\newcommand{\tM}{\operatorname{\widetilde{\mathcal{M}}}}
\newcommand{\tLoc}{{\widetilde{Loc}}}
\newcommand{\Ll}{\operatorname{{\mathcal{L}_\lambda}}}
\newcommand{\tL}{\operatorname{\widetilde{\mathcal{L}}}}
\newcommand{\tU}{{\widetilde{U}_q}}
\newcommand{\tTX}{\operatorname{\widetilde{T}^\star X}}
\newcommand{\Ind}{\operatorname{Ind}}

\newcommand{\Hom}{\operatorname{Hom}}
\newcommand{\End}{\operatorname{End}}

\newcommand{\gr}{\operatorname{gr}}

\newcommand{\mood}{\hbox{\ensuremath{
\operatorname{mod}}}}

\newcommand{\ad}{\hbox{\ensuremath{\operatorname{ad}}}}
\newcommand{\mO}{\mathcal{O}}

\newcommand{\Z}{\ensuremath{\mathcal{Z}}}
\newcommand{\Zl}{\ensuremath{\mathcal{Z}^{(l)}}}
\newcommand{\ZHC}{\ensuremath{\mathcal{Z}^{HC}}}
\newcommand{\Znl}{\ensuremath{\mathcal{Z}^{(l)}_{-,\lambda}}}
\newcommand{\Zbl}{\ensuremath{\mathcal{Z}^{(l)}_{\mathfrak{b}_-}}}

\newcommand{\C}{\ensuremath{{\mathbb{C}}}}
\newcommand{\D}{\ensuremath{\mathcal{D}_{G_q}}}
\newcommand{\Dl}{\ensuremath{\mathcal{D}^\lambda_q}}

\newcommand{\g}{\ensuremath{\mathfrak{g}}}

\newcommand{\n}{\ensuremath{\mathfrak{n}}}
\newcommand{\h}{\ensuremath{\mathfrak{h}}}

\newcommand{\bb}{\ensuremath{\mathfrak{b}}}

\newcommand{\Oq}{\ensuremath{\mathcal{O}_q}}

\newcommand{\UA}{\ensuremath{U_\mathcal{A}}}
\newcommand{\UAres}{\ensuremath{U^{res}_\mathcal{A}}}

\newcommand{\Fp}{\ensuremath{\mathbb{F}_p}}
\newcommand{\UBr}{\ensuremath{U^{res}_q(\bb)}}
\newcommand{\W}{\ensuremath{\mathcal{W}}}
\newcommand{\Tl}{\ensuremath{T^\star X^{\lambda}}}
\newcommand{\M}{\ensuremath{\mathcal{M}}}
\newcommand{\MBG}{\ensuremath{\mathcal{M}_{B_q}(G_q)}}

\newcommand{\DBG}{\ensuremath{\mathcal{D}^\lambda_{B_q}(G_q)}}
\newcommand{\tDBG}{\ensuremath{\widetilde{\mathcal{D}}_{B_q}(G_q)}}

\newcommand{\U}{\ensuremath{U^{fin}_q}}
\newcommand{\Ul}{\ensuremath{U^{\lambda}_q}}

\newcommand{\A}{\ensuremath{\mathcal{A}}}

\newcommand{\OA}{\ensuremath{\mathcal{O}_\mathcal{A}}}

\theoremstyle{plain}

\newtheorem{Thm}{Theorem}[section]
\newtheorem{Prop}[Thm]{Proposition}
\newtheorem{Lem}[Thm]{Lemma}

\newtheorem{Cor}[Thm]{Corollary}

\newtheorem{DT}[Thm]{Definition-Proposition}

\theoremstyle{definition}
\newtheorem{defi}[Thm]{Definition}

\theoremstyle{remark}
\newtheorem{Rem}[Thm]{Remark}

\begin{document}
\title{Localization for quantum groups at a root of unity.}

\author{Erik Backelin and Kobi Kremnitzer}
\maketitle

\begin{abstract} In the paper \cite{BK} we defined categories of
equivariant quantum $\mathcal{O}_q$-modules and
$\mathcal{D}_q$-modules on the quantum flag variety of $G$. We
proved that the Beilinson-Bernstein localization theorem holds at
a generic $q$. Here we prove that a derived version of this
theorem holds at the root of unity case. Namely, the global
section functor gives a derived equivalence between categories of
$U_q$-modules and $\mathcal{D}_q$-modules on the quantum flag
variety.

For this we first prove that $\mathcal{D}_q$ is an Azumaya algebra
over a dense subset of the cotangent bundle $T^\star X$ of the classical (char $0$)
flag variety $X$. This way we get a derived equivalence between
representations of $U_q$ and certain $\mathcal{O}_{T^\star
X}$-modules.

In the paper \cite{BMR} similar results were obtained for a Lie
algebra $\g_p$ in char $p$. Hence, representations of $\g_p$ and
of $U_q$ (when $q$ is a p'th root of unity) are related via the
cotangent bundles $T^\star X$ in char $0$ and in char $p$,
respectively.
\end{abstract}

\section{Introduction}
Let $\C$ be the field of complex numbers and fix $q \in \C^\star$.
Let $\g$ be a semi-simple Lie algebra over $\C$ and let $G$ be the
corresponding simply connected algebraic group. Let $U_q$ be a
quantized enveloping algebra of $\g$. Let $\Oq$ be the algebra of
quantized functions on $G$. Let $\Oq(B)$ be the quotient Hopf
algebra of $\Oq$ corresponding to a Borel subgroup $B$ of $G$.

\smallskip

\noindent \textbf{We now recall the main results in \cite{BK}}.
The constructions given there are crucial for the present paper
and a fairly detailed survey of the material there is given in the
next section. We defined an equivariant sheaf of quasicoherent
modules over the quantum flag variety to be a left $\Oq$-module
equipped with a right $\Oq(B)$-comodule structure satisfying
certain compatibility conditions. Such objects form a category
denoted $\MBG$. It contains certain line bundles $\Oq(\lambda)$
for $\lambda$ in the weight lattice. We proved that $\Oq(\lambda)$
is ample for $\lambda >> 0$ holds for every $q$. This implies that
the category $\MBG$ is a Proj-category in the sense of Serre.

We defined the category $\DBG$ of $\lambda$-twisted quantum
$\mathcal{D}$-modules on the quantized flag variety (see
definition \ref{d12}). An object $M \in \DBG$ is an object in
$\MBG$ equipped with an additional left $U_q$-action satisfying
certain compatibility axioms (in particular, the
$U_q(\mathfrak{b})$-action on $M$, that is the restriction of the
$U_q$-action, and the $\Oq(B)$-coaction "differ by $\lambda$".

The global section functor $\Gamma$ on $\MBG$ (and on $\DBG$) is
given by taking $\Oq(B)$-coinvariants. Let $M \in \MBG$. Its
sections over an open set (i.e., a localization $\mathcal{O}_{qf}$
of $\Oq$) are given by $\Gamma(\mathcal{O}_{qf}, M) =
\Gamma(\mathcal{O}_{qf} \otimes_{\Oq} M)$.

The category $\DBG$ has a distinguished object $\Dl$. As a vector
space $\Dl = \Oq \otimes M_\lambda$, where $M_\lambda$ is a Verma
module. $\Dl$ is the (quantum) equivariant counterpart of the
usual sheaf of rings of $\lambda$-twisted differential operators
on the flag variety of $G$. (The "set" $\Dl$ is not a ring, but
its "sections over open sets" are naturally rings.) We proved that
for each $q$ except a finite set of roots of unity (depending on
$\g$), the global sections $\Gamma(\Dl)$ is isomorphic to $\Ul :=
\U/J_\lambda$, where $J_\lambda$ is the annihilator of $M_\lambda$
and $U^{fin}_q$ is the maximal subalgebra of $U_q$ on which the
adjoint action of $U_q$ is locally finite.

The main result in \cite{BK} stated that the global section
functor gives an equivalence of categories between $\DBG$ and
modules over $\Ul$ in the case when $q$ is not a root of unity and
$\lambda$ is regular. This is a quantum version of
Beilinson-Bernsteins, \cite{BB}, localization theorem.

\smallskip
\noindent \textbf{In the present paper} we study the root of unity
case, $q$ is a primitive $l$'th root of unity. The situation
becomes different and very interesting: It turns out that $\Dl$
naturally forms a sheaf of algebras  over the classical
(non-quantum) complex variety $T^\star X^\lambda := (G\times
N^{\lambda^{2l}})/B$, where $N^{\lambda^{2l}}$ is the
$B$-submodule $\{\lambda^{2l}\} \times N$ of $B$ (see section
$3.1.2$). If $\lambda$ is integral $N^{\lambda^{2l}} = N$ and so
$T^\star X^\lambda = T^\star X$ is the cotangent bundle of the
flag variety $X$ of $G$ in this case. Hence, we refer to $T^\star
X^\lambda$ as a twisted cotangent bundle.

A key observation is that $\Dl$ is an Azumaya algebra over a dense subset of $\Tl$
(proposition \ref{Azumaya}) and that this Azumaya algebras splits
over formal neighborhoods of generalized Springer fibers
(proposition \ref{Azsplit}). We then show that that $\Dl$ has no
higher self extensions, i.e. that $R\Gamma^{>0}(\Dl) = 0$
(proposition \ref{RGammaD}). This fact, together with the result
$\Gamma(\Dl) = \Ul$ and the Azumaya property implies a derived
version of the localization theorem: The functor $\Gamma$ induces
an equivalence between bounded derived categories $D^b(\DBG)$ and
$D^b(\Ul -\mood)$ (restricted to the dense subsets).

Using the Azumaya splitting we prove that the subcategory of
$\DBG$ whose objects are supported on (a formal neighbourhood of a
generalized) Springer fiber is isomorphic to the category of
$\mO$- modules over the twisted cotangent bundles $T^\star
X^\lambda$ supported on the same fiber (corollary \ref{AZTX}).

Combining these results we get an equivalence between $D^b(\Ul
-\mood)$ and the derived category of $\mO$-modules on $\Tl$ whose
cohomologies are supported on certain Springer fibers.

An application of our theory is for instance the computation of
the number of simple $\mathfrak{u}_q$-modules, because this number
can be interpreted as the rank of the $K$-group of the category of
$\mO$-modules on $\Tl$ supported on the trivial Springer fiber.
Of course, such a formula follows also from the link to the
representation theory of an affine Kac-Moody Lie algebras,
established by Kazhdan and Lusztig. Our method gives the possibility to extend this to nontrivial central characters as well. In a future paper we will use this to prove a conjecture of DeConcini, Kac and Processi regarding dimensions of irreducible modules.
\smallskip

\noindent \textbf{In \cite{BMR} analogous results were established
for a Lie algebra $\g_p$ in characteristic $p$. } In fact most of
our methods of proofs are borrowed from that paper. They showed
analogous results for $U(\g_p)$-modules, (Crystalline)
$\mathcal{D}$-modules on the flagvariety $X(\bar{\mathbb{F}}_p)$
and certain twisted cotangent bundles of $X$ over
$\bar{\mathbb{F}}_p$.

Combining their results with ours we see that the representation
theory of $U_q$ (when $q$ is a p'th root of unity) is related to
the representation theory of $U(\g_p)$ via cotangent bundles of
$X$ in $char \, 0$ and $char \, p$, respectively. We know
furthermore that baby-Verma modules go to Skyscraper sheaves in
both cases and \cite{BMR} showed that the $K$-groups of the
cotangent bundle categories are isomorphic. Right now we are
investigating what can be deduced about the representation theory
of $\g_p$, e.g. character formulas, from the representation theory
of $U_q$ with these methods.

\subsection{Acknowledgements} We thank Joseph Bernstein
and David Kazhdan for many useful conversations.
\section{Generalities}

\subsection{Quantum groups} See Chari and Pressley \cite{CP} for details
about the topics in this section.
\subsubsection{Conventions}\label{Conventions}

Let $\C$ be the field of complex numbers and fix $q \in \C^\star$.

\textbf{We always assume that if $q$ is a root of unity it is
primitive of odd order and in case $G$ has a component of type
$G_2$ the order is also prime to 3.} Let $\A$ be the local ring
$\mathbb{Z}[\nu]_{\mathfrak{m}}$, where $\mathfrak{m}$ is the
maximal ideal in $\mathbb{Z}[\nu]$ generated by $\nu - 1$ and a
fixed odd prime $p$.

\subsubsection{Root data} Let $\g$ be a semi-simple Lie
algebra and let $\h \subset \bb$ be a Cartan subalgebra contained
in a Borel subalgebra of $\g$. Let $R$ be the root system, $\Delta
\subset R_+ \subset R$ a basis and the positive roots. Let $P
\subset \h^\star$ be the weight lattice and $P_+$ the positive
weights; the $i$'th fundamental weight is denoted by $\omega_i$
and $\rho$ denotes the half sum of the positive roots. Let $Q
\subset P$ be the root lattice and $Q_+ \subset Q$ those elements
which have non-negative coefficients with respect to the basis of
simple roots. Let $\mathcal{W}$ be the Weyl group of $\g$. We let
$<\,,\,>$ denote a $\mathcal{W}$-invariant bilinear form on
$\h^\star$ normalized by $<\gamma,\gamma> = 2$ for each short root
$\gamma$.

Let $T_P = \Hom_{groups}(P,k^\star)$ be the character group of $P$
with values in $k$ (we use additive notation for this group). If
$\mu \in P$, then $<\mu,P> \subset \mathbb{Z}$ and hence we can
define $q^\mu \in T_P$ by the formula $q^\mu(\gamma) =
q^{<\mu,\gamma>}$, for $\gamma \in P$. If $\mu \in P, \lambda \in
T_P$ we write $\mu + \lambda = q^\mu + \lambda$. Note that the
Weyl group naturally acts on $T_P$.

\subsubsection{Quantized enveloping algebra $U_q$ and quantized
algebra of functions $\Oq$.} Let $U_q$ be the simply connected
quantized enveloping algebra of $\g$ over $\C$. Recall that $U_q$
has algebra generators $E_\alpha, F_\alpha, K_\mu$, $\alpha,
\beta$ are simple roots, $\mu \in P$ subject to the relations
\begin{equation}\label{qgr1}
K_{\lambda} K_{\mu} = K_{\lambda+\mu},\;\; K_0 = 1,
\end{equation}
\begin{equation}\label{qgr2}K_\mu E_\alpha K_{-\mu} = q^{<\mu, \alpha>} E_\alpha,\;\; K_\mu F_\alpha K_{-\mu} = q^{-<\mu, \alpha>}
F_\alpha,\end{equation}
\begin{equation}\label{qgr3}[E_\alpha,\, F_\beta] = \delta_{\alpha, \beta} {{K_\alpha - K_{-\alpha}}\over{q_\alpha-q^{-1}_\alpha}}\end{equation}
and certain Serre-relations that we do not recall here. Here

$q_\alpha = q^{d_\alpha}$, $d_\alpha = <\alpha,\alpha>/2$. (We
have assumed that $q^2_\alpha \neq 1$.)

Let $G$ be the simply connected algebraic group with Lie algebra
$\g$, $B$ be a Borel subgroup of $G$ and $N \subset B$ its
unipotent radical. Let $\bb = \operatorname{Lie} B$ and
$\mathfrak{n} = \operatorname{Lie} N$ and denote by $U_q(\bb)$ and
$U_q(\mathfrak{n})$ the corresponding subalgebras of $U_q$. Then
$U_q(\bb)$ is a Hopf algebra, while $U_q(\mathfrak{n})$ is only an
algebra. Let $\Oq = \Oq(G)$ be the algebra of matrix coefficients
of finite dimensional type-1 representations of $U_q$. There is a
natural pairing $(\;,\;): U_q \otimes \Oq \to \C$. This gives a
$U_q$-bimodule structure on $\Oq$ as follows
\begin{equation}\label{23}
ua = a_1 (u,a_2),\;\; au = (u,a_1)a_2, \,\, u \in U_q, a \in \Oq
\end{equation}
Then $\Oq$ is the (restricted) dual of $U_q$ with respect to this
pairing. We let $\Oq(B)$ and $\Oq(N)$ be the quotient algebras of
$\Oq$ corresponding to the subalgebras $U_q(\bb)$ and
$U_q(\mathfrak{n})$ of $U_q$, respectively, by means of this
duality. Then $\Oq(B)$ is a Hopf algebra and $\Oq(N)$ is only an
algebra.

There is a braid group action on $U_q$. For each $w \in \W$, we
get an automorphism $T_w$ of $U_q$.

\subsubsection{Integral versions of $U_q$.} Let
$\UAres$ be the Lusztig's integral form of $U_q$, the $\A$-algebra
in $U_q$ generated by divided powers $E^{(n)}_\alpha =
E^n_\alpha/{[n]_{d_\alpha}!}$, $F^{(n)}_\alpha =
F^n_\alpha/{[n]_{d_\alpha}!}$, $\alpha$ a simple root, $n \geq 1$
(where $[m]_d = {\prod^m_{s=1} {q^{d\cdot s} -  q^{-d\cdot
s}}\over{q^{d} -  q^{-d}}}$) and the $K_\mu$'s, $\mu \in P$. There
is also the De Consini-Kac integral form $\UA$, which is generated
over $\A$ by the $E_\alpha, F_\alpha$ and $K_\mu$'s. The
subalgebra $\UA$ is preserved by the adjoint action of $\UAres$:
$ad_{\UAres} (\UA) \subset \UA$. The operators $T_w$ from section
2.1.2 preserves the integral versions.

$\OA$ is defined to be the dual of $\UAres$. This is an $\A$-sub
Hopf algebra of $\Oq$.

\subsubsection{Finite part of $U_q$.} The algebra $U_q$ acts on
itself by the adjoint action $\ad: U_q \to U_q$ where $\ad(u)(v) =
u_1 v S(u_2)$. Let $U^{fin}_q$ be the finite part of $U_q$ with
respect to this action:
$$ U^{fin}_q = \{v \in U_q; \dim \ad(U_q)(v) < \infty\}.$$
This is a subalgebra. (See \cite{JL}.)

We can also give an integral version of the finite part as the
finite part of the action of $\UAres$ on $\UA$.Thus by
specializing we get a subalgebra of $U_q$ for every $q$. Of
course, when specialized to generic $q$ this coincides with the
previous definition.

\subsubsection{Specializations and Frobenius maps.}
For any ring map $\phi: \A \to R$ we put $U_R = \UA \otimes_\A R$
and $U^{res}_R = \UAres \otimes_\A R$. If $R = \C$ and $\phi(\nu)
= q$, there are three different cases: $q$ is a root of unity, $q
= 1$ and $q$ is generic. Then $U_R = U_q$.

There is the also the ring map $\A \to \Fp$, sending $\nu \to 1$.
Then $U_{\Fp} /(K-1) = U(\g_p)$, the enveloping algebra of the Lie
algebra $\g_p$ in characteristic $p$.

For any $\UA$-module (resp. $\UAres$-module) $M_\A$ we put $M_R =
M_\A \otimes_\A R$. This is an $U_R$-module (resp.
$U^{res}_R$-module). When $R = \C$ we simply write $M = M_\C$.

When $q$ is a root of unity, we have the Frobenius map: $U^{res}_q
\to U(\g)$. Its algebra kernel is denoted by $\mathfrak{u}_q$. We
also have the Frobenius map $U^{res}_q(\bb) \to U(\bb)$, with
algebra kernel $b_q$. These maps induces dual maps $\mO = \mO(G)
\hookrightarrow \Oq$ and $\mO(B) \hookrightarrow \Oq(B)$.

For each $q$ there exists a map $U_q \to U^{res}_q$ whose image is
$\mathfrak{u}_q$ and whose algebra kernel is $\Zl$ (see section
2.1.7 below for the definition of $\Zl$).

\subsubsection{Verma modules.}\label{Verma} For each $\lambda \in T_P$ there is the
one dimensional $U_q(\bb)$-module $\C_\lambda$ which is given by
extending $\lambda$ to act by zero on the $E_\alpha$'s. The
Verma-module $M_\lambda$ is the $U_q$-module induced from
$\C_\lambda$. If $\mu \in P$ we write $M_\mu = M_{q^\mu}$. A point
important for us is that $M_\lambda$ carries an
$U^{res}_q(\bb)$-module structure defined as follows:
$U^{res}_q(\bb)$ acts on $U_q$ by restricting the adjoint action
of $U^{res}_q$ on $U_q$. This induces a $U_q(\bb)$-action on the
quotient $M_\lambda$ of $U_q$. Since this action is locally finite
it corresponds to an $\Oq(B)$-comodule action on $M_\lambda$.
\textbf{Note!!} As a $U^{res}_q(\bb)$-module $M_\lambda$ has
\textit{trivial highest weight}. (In case $q$ is generic,
$U^{res}_q(\bb) = U_q(\bb)$ and then the $U^{res}_q(\bb)$ action
on $M_\lambda$ described above is the same as the $U(\bb)$-action
on $M_\lambda \otimes \C_{-\lambda}$.)

Verma module $M_\lambda$ has an integral version $M_{\lambda,\A}$.

\subsubsection{Centers of $U_q$ and definition of $\tU$.} Let $\Z$ denote the center of $U_q$.
When $q$ is a $p$'th root of unity $\Z$ contains the
Harish-Chandra center $\ZHC$ and the $l$'th center $\Zl$ which is
generated by the $E^{l}_\alpha$, $F^{l}_\alpha$, $K^{l}_\mu$ and
$K^{-l}_\mu$'s. In fact, $\Z = \Zl \otimes_{\Zl \cap \ZHC} \ZHC$.
There is the Harish-Chandra homomorphism $\ZHC \to \mO(T_P)$ that
maps isomorphically to the $W$-invariant even part of
$\mO(T_{P})$. We define $\tU = U_q \otimes_{\ZHC} \mO(T_P)$.

\subsubsection{Some conventions.} We shall frequently refer to a right
(resp. left) $\Oq$-comodule as a left (resp. right) $G_q$-module,
etc. If we have two right $\Oq$-comodules $V$ and $W$, then
$V\otimes W$ carries the structure of a right $\Oq$-comodule via
the formula
$$
\delta(v\otimes w) = v_1\otimes w_1 \otimes v_2w_2$$ We shall
refer to this action as the {\it tensor} or {\it diagonal} action.
A similar formula exist for left comodules.

\subsection{Quantum flag variety}
Here we recall the definition and basic properties of the quantum
flag variety from \cite{BK}.

\subsubsection{Category $\MBG$.} The composition
\begin{equation}\label{e2}
\Oq \to \Oq \otimes \Oq \to \Oq \otimes \Oq(B)\end{equation}
defines a right $\Oq(B)$-comodule structure on  $\Oq$. A
$B_q$-equivariant sheaves on $G_q$ is a triple $(F, \alpha,
\beta)$ where $F$ is a vector space, $\alpha: \Oq \otimes F \to F$
a left $\Oq$-module action and $\beta: F \to F \otimes \Oq(B)$ a
right $\Oq(B)$-comodule action such that $\alpha$ is a right
comodule map, where we consider the tensor comodule structure on
$\Oq(G) \otimes F$.
\begin{defi}\label{hejsan} We denote $\MBG$ to be the category of $B_q$-equivariant
sheaves on $G_q$. Morphisms in $\MBG$ are those compatible with
all structures.\end{defi}

If $q = 1$, the category $\M_B(G)$ is equivalent to the category
$\M(G/B)$ of quasi-coherent sheaves on $G/B$.

\begin{defi} We define the induction functor $\Ind: B_q-\mood$ to $\MBG$, $\Ind
V = \Oq \otimes V$ with the tensor $B_q$-action and the
$\Oq$-action on the first factor. For $\lambda \in P$ we define a
line bundle $\Oq(\lambda) = \Ind \C_{-\lambda}$.
\end{defi}
\begin{defi} The global section functor $\Gamma: \MBG \to \C-\mood$
is defined by
$$\Gamma(M) = \Hom_{\MBG}(\Oq, M) = \{m \in M; \Delta_B(m) =
m\otimes 1\}.$$ This is the set of $B_q$-invariants in $M$.
\end{defi}

The category $\MBG$ has enough injectives, so derived functors
$R\Gamma$ are well-defined. We showed that $R^i\Gamma(\Ind V) =
H^i(G_q/B_q,V)$, where $H^i(G_q/B_q,\;)$ is the $i$'th derived
functor of the functor $V \to \Gamma(\Ind V)$ from $B_q-\mood$ to
$\C-\mood$.

We proved a quantum version of Serre's basic theorem on projective
schemes: Each $M \in \MBG$ is a quotient of a direct sum of
$\Oq(\lambda)$'s and each surjection $M \twoheadrightarrow M'$ of
noetherian objects in $\MBG$ induces a surjection
$\Gamma(M(\lambda)) \twoheadrightarrow \Gamma(M'(\lambda))$ for
$\lambda >> 0$.

Here the notation $\lambda >> 0$ means that $<\lambda,
\alpha^\wedge>$ is a sufficiently large integer for each simple
root $\alpha$ and $M(\lambda) = M \otimes \C_{-\lambda}$ is the
$\lambda$-twist of $M$.

Let $V \in G_q-\mood$. Denote by $V \vert B_q$ the restriction of
$V$ to $B_q$ and by $V^{triv}$ the trivial $B_q$-module whose
underlying space is $V$. We showed that $\Ind V \vert B_q$ and
$\Ind V^{triv}$ are isomorphic in $\MBG$. In particular
\begin{equation}\label{Gmod}
\Gamma(\Ind V \vert B_q) = V \vert B_q \otimes \Gamma(\Oq) = V
\vert B_q, \hbox{ for } V \in G_q-\mood
\end{equation}
\subsubsection{ \MBG at a root of unity}
In case $q$ is a root of unity we have the following Frobenius
morphism:
\begin{equation}
Fr_*:\MBG \rightarrow \M(G/B)
\end{equation}
defined as
\begin{equation}
 N \mapsto N^{b_q}
\end{equation}

Using the description of $\MBG$ as $Proj(A_q)$ where $A_q=
\bigoplus V_{q,\lambda}$ \cite{BK} and similarly $\M(G/B)=
Proj(A)$ where $A= \bigoplus V_\lambda$, we see that $Fr_*$ is
induced from the quantum Frobenius map $A\hookrightarrow A_q$. It
follows that:
\begin{Prop}\label{Frob} $Fr_*$ is exact and faithful.
\end{Prop}

 This functor has a left adjoint.

\section{$\Dl$-modules at a root of unity }
  \subsection{ First construction}\label{actionofG}
From now on $q$ is an $l$'th root of unity (recall the
restrictions of \ref{Conventions}).

 In this section we shall give a representation theoretic
 construction of the sheaf of quantum differential operators. This
 will turn out to be a sheaf of algebras over the Springer
 resolution- the sheaf of endomorphisms of the (nonexistent)
 universal baby Verma module.

 Recall the Frobenius map $O
\hookrightarrow O_q$ and the fact that $O_q^{\mathfrak{u}_q}=O$.
This allows us to define the functor of (finite) induction
\begin{equation}
\Ind:\mathfrak{u}_q-mod \rightarrow
\ensuremath{\mathcal{M}_{G_q}(G)}
\end{equation}
\begin{equation}
 N\mapsto (O_q\otimes N)^{\mathfrak{u}_q}
\end{equation}

Here  \ensuremath{\mathcal{M}_{G_q}(G)} is the category of $G_q$
equivariant $O$-modules, that is: an $O$ module which is also an
$O_q$ comodule, and such that the $O$ module structure map is a
map of $O_q$ comodules.

We have the following important proposition \cite{AG}:
\begin{Prop}\label{smalluq}
 $Ind:\mathfrak{u}_q-mod \rightarrow
\ensuremath{\mathcal{M}_{G_q}(G)}$ is an equivelence of
categories.
\end{Prop}

Notice that both categories are tensor categories (in
$\ensuremath{\mathcal{M}_{G_q}(G)}$ it is tensoring over $O$) and
that $Ind$ is a tensor functor. This will be used later. Note also
that $\ensuremath{\mathcal{M}_{G_q}(G)}$ has an obvious action by
$G$, that is for any $g \in G$ we have an functor
$F_g:\ensuremath{\mathcal{M}_{G_q}(G)}\rightarrow
\ensuremath{\mathcal{M}_{G_q}(G)}$ and natural transformations
$\alpha_{g,h}:F_g \circ F_h \Rightarrow F_{gh}$ satisfying a
certain cocycle condition. Hence starting from any
$\mathfrak{u}_q$ module we can form a family of such modules
indexed by $G$, more precisely an $O(G)$-module in the category of
$\mathfrak{u}_q$- modules. This will give the `universal family`
of baby verma modules.

 \cite{CKP} defined an action of an infinite dimensional group $\mathfrak{G}$ on $U_q$ 
preserving the l-center (and the augmetation ideal of the l-center) and thus acting also on $\mathfrak{u}_q$. This is defined by observing that
the derivation defined by commuting with the divided powers $E^{(l)}, F^{(l)}$ actually preserves the algebra generated by the nondivided powers. These derivations are then exponentiated to get automorphisms of $U_q$ at a root of unity. The group they generate is infinite dimensional as the action is not locally finite. This action also induces an action on $\mathfrak{u}_q$ but here the group is finite dimensional $\mathfrak{G}_0$ . We thus have another group action on the category.   

\begin{Prop} \label{sameaction} Let $r:Aut(\mathfrak{u}_q)\to Out(\mathfrak{u}_q)$ be the natural map. $r(\mathfrak{G}_0)=G$. 
\end{Prop}

\textit{proof of proposition \ref{sameaction}} It is enough to check that the action on the category is the same. Since both actions are on a category of modules over a finite dimensional algebra it is enough to check that they agree infinitessimely. That is: any action defines a map from the Lie algebra to outer derivations of the algebra and it is enough to check on this level. But for  $\mathfrak{G}_0$  the Lie algebra action is given by the derivation defined by commuting with a divided power and for $G$ it is given by the adjoint action and both have the same image inside $Ext^1(\mathfrak{u}_q,\mathfrak{u}_q)$ (they both span outer derivations).

Recall that $M_\lambda = U_q / J_\lambda$; put
\begin{equation}
I_\lambda = \Zl \cap J_\lambda
\end{equation}
\begin{equation}
 \Znl = \Zl/ I_\lambda
\end{equation}
 so that
$\Znl \subset M_\lambda$.

 Then $\Znl$ is a $U^{res}_q(\bb)$-module
algebra and $M_\lambda$ is a $U^{res}_q(\bb)$-module for this
algebra (see section \ref{Verma}).

In fact, $\mO \otimes \Znl$ is a $B$-equivariant algebra: Recall
the Frobenius map $Fr: U^{res}_q(\bb) \to U(\bb)$ and denote the
algebra kernel of  $Fr$ by $b_q$.

\begin{defi}\label{twists} $i)$ Consider $\mO(B)$ as a $B$-module under the adjoint action.
For each $t \in T$ we have the $B$-submodule $N^t = \{t\} \times
N$ of $B$. $\mO(N^t)$ is isomorphic to $\mO(N)$ as an algebra, but
not as a $B$-module, unless $t = 1$.

$ii)$ Think of $T$ as characters on $2l \cdot P$ and of $T_P$ as
characters on $P$. Lattice inclusion induces a natural map
${(\;)}^{2l}: T_P \to T$.
\end{defi}
Note that for $\lambda$ integral $\lambda^{2l} = 1$.

 We now have
\begin{Lem}\label{TXlem}  $b_q$ acts trivially on
$\Znl$, this module is a pullback by the Frobenius of the
$B$-module $\mO(N^{\lambda^{2l}})$.
\end{Lem}
\smallskip

\noindent \textit{Proof of proposition \ref{TXlem}.} This will
follow from $\ref{Bmod}$. $\Box$
\smallskip

 $b_q$ also acts trivially on $\mO$, so $\mO
\otimes \Znl$ is a $B$-module. Hence, we have the category
$\mathcal{M}_B(\mO \otimes \Znl)$ of $B$-equivariant $\mO \otimes
\Znl$-modules; from lemma \ref{TXlem} we conclude that

\begin{equation}\label{equiv}
\mathcal{M}_B(\mO \otimes \Znl) \cong qcoh(\mO_{(G\times
N^{\lambda^{2l}})/B})
\end{equation}

 Here the $B$ action on $G\times N^\lambda$ is given by $b \cdot
(g, x) = (bg, b\cdot x)$. Similarly we can define the category
${}_{G_q}\mathcal{M}_B(\mO \otimes \Znl)$ of $G_q$-equivariant
objects in $\mathcal{M}_B(\mO \otimes \Znl)$.

We shall denote
\begin{equation}
G\times N^{\lambda^{2l}}/B = \Tl
\end{equation}

Note that for integral $\lambda$ this is the cotangent bundle to
the flag variety, also known as the Springer resolution. For non
integral $\lambda$ this is a twisted cotangent bundle.

For any $\tau \in maxspec(\Znl)$ we have the central reduction
$M_{\lambda,\tau}$ (a baby Verma module). Only for trivial $\tau$
(corresponding to the augmentation ideal) we get a
$\mathfrak{u}_q$-module. But for any $\tau$,
$End(M_{\lambda,\tau})=M_{\lambda,\tau}\otimes M_{\lambda,\tau}^*$
is a $\mathfrak{u}_q$-module since its l-central character is
trivial, and likewise $End_{\Znl} (M_\lambda)$ is a
$\mathfrak{u}_q$-module. Hence we can define:

\begin{defi}\label{firstdefi}
 $D = (Ind(End_{\Znl} (M_\lambda+2\rho)))$.
\end{defi}
(The shift by $2\rho$ will become clear later.)
Since $End_{\Znl} (M_\lambda)$ is a $\Znl$-module and also a
$B_q$-module ($O_q(B)$-comodule) in a compatible way we get that
$D$ actually lives in ${}_{G_q}\mathcal{M}_B(\mO \otimes \Znl)$.
Since $Ind$ is a tensor functor and $End_{\Znl} (M_\lambda)$ is a
$\mathfrak{u}_q$ algebra we get that $D$ is a sheaf of algebras
over the (twisted) cotangent bundle.

By construction we know that the algebras sitting over the fiber
over $B$ are matrix algebras (endomorphisms of baby Verma
modules), hence we get
\begin{Prop}\label{Azumaya}
 Over a dense subset of $\Tl$, $\Dl$ is an
Azumaya algebra.
\end{Prop}

Note that this dense subset contains the zero section (and the Azumaya algebra is even split over it) since for modules with trivial l-central character we have an action already on them and not only on their endomorphisms.

\begin{Rem}
Note that all our constructions can also be defined over a formal neighbourhood of a prime $p$ that is over a p-adic field and that when specialized to $\mathbf{F}_p$ they give the usual characteristic p crystalline differential opeartors which are Azumaya and thus we would get that over the p-adic field our algebra is Azumaya as well. This will not be used in this paper since we will look at complex representations, but in a subsequent paper we would use this to construct t-structures in zero characteristic relating to the ones constructed by Bezrukavnikov, Mircovic and Rumynin \cite{BMR}. 
\end{Rem}

\begin{Rem}
For any rigid braided tensor category one can define the notion of
an Azumaya algebra. In the category of $\mathfrak{u}_q$-modules
$End_{\Znl} (M_\lambda)$ is an Azumaya algebra. Hence, using the
equivalence \ref{smalluq} we get that $D$ is an Azumaya algebra
over $\Tl$, not with respect to the usual braiding (flip) but with
respect to the braiding induced from $U^{res}_q$.
\end{Rem}

 \subsection{Second construction-The ring $\D$ and the category of $\Dl$
modules} We need the following important
\begin{Rem}\label{Integral versions.} All objects described in the preceding chapters are
defined over $\A$. For any specialization $\A \to R$ and any
object $Obj$ we denote by $Obj_R$ its $R$-form. For the functors
we don't use any subscripts; so, for instance, there is the
functor $\Ind: B_{q,R}-\operatorname{mod} \to \MBG_R$.
\end{Rem}
Recall the $U_q$-bimodule structure on $\Oq$ given by \ref{23}.
Now, as we have two versions of the quantum group we pick the
following definition of the ring of differential operators on the
group (the $crystalline$ version).
\begin{defi} We define the ring of quantum differential operators on $G_q$ to be
the smash product algebra $\D := \Oq \star U_q$. So $\D = \Oq
\otimes U_q$ as a vector space and multiplication is given by
\begin{equation}\label{e71}
a \otimes u \cdot b \otimes v = au_1(b)\otimes u_2v.
\end{equation}
\end{defi}
We consider now the ring $\D$ as a left $U^{res}_q$-module, via
the left $U^{res}_q$-action on $\Oq$ in \ref{23} and the left
adjoint action of $U^{res}_q$ on $U_q$; this way $\D$ becomes a
module algebra for $U^{res}_q$: In the following we will use the
restriction of this action to $U^{res}_q(\bb) \subset U^{res}_q$.
As $U_q$ is not locally finite with respect to the adjoint action,
this $U^{res}_q(\bb)$-action doesn't integrate to a $B_q$-action.
Thus $\D$ is not an object of $\MBG$; however, $\D$ has a
subalgebra $\D^{fin} = \Oq \star \U$ which belongs to $\MBG$. This
fact will be used below.

\begin{defi}\label{d12} Let $\lambda \in T_P$. A $(B_{q},\lambda)$-equivariant $\D$-module
is a triple $(M, \alpha, \beta)$, where $M$ is a $\C$-module,
$\alpha: \D \otimes M \to M$ a left $\D$-action and $\beta^{res}:
M \to M\otimes \Oq(B)$ a right $\Oq(B)$-coaction. The latter
action induces an $U^{res}_q(\bb)$-action on $M$ again denoted by
$\beta^{res}$. We have the natural map $U_q(\bb) \to \UBr$ which
together with $\beta^{res}$ gives an action $\beta$ of $U_q(\bb)$
on $M$. We require
\smallskip

\noindent $i)$ The $U_q(\bb)$-actions on $M \otimes \C_\lambda$
given by $\beta \otimes \lambda$ and by $(\alpha\vert_{U_q(\bb)})
\otimes \operatorname{Id}$ coincide.

\noindent $ii)$ The map $\alpha$ is $U^{res}_q(\bb)$-linear with
respect to the $\beta$-action  on $M$ and the action on $\D$.

\smallskip

These objects form a category denoted $\DBG$. There is the
forgetful functor $\DBG \to \MBG$. Morphisms in $\DBG$ are
morphisms in $\MBG$ that are $\D$-linear.
\end{defi}
We defined $\Dl$ as the maximal quotient of $\D$ which is an
object of $\DBG$ and showed that
\begin{equation}\label{relevantstruct1}
\Dl = \Ind M_\lambda
\end{equation}
as an object in $\DBG$. (See section \ref{Verma} for the
$B_q$-action $=$ $U^{res}_q(\bb)$-action on $M_\lambda$.). The
global section functor $\Gamma: \DBG \to \M$ is the functor of
taking $B_q$ invariants (with respect to the action $\beta$); we
have $\Gamma = \Hom_{\DBG}(\Dl,\;)$.

Hence, in particular $\Gamma(\Dl) = \End_{\DBG}(\Dl)$ (which
explains the ring structure on $\Gamma(\Dl)$).

If we view $\Dl$ as an object in ${}_{G_q}\MBG$ we have that
$Fr_*(\Dl)$ is an object of ${}_{G_q}\M(G/B)$. But it actually
lies in ${}_{G_q}\mathcal{M}_B(\mO \otimes \Znl)$ since
$M_\lambda$ is a $\Znl$ module. In the next section we shall show
that the two constructions of $\Dl$ coincide and that the category
of modules over this sheaf of algebras is $\DBG$.

\subsection{$\Dl$ as a sheaf of algebras}
\smallskip
We have a natural functor:
\begin{equation}
F: \DBG \rightarrow \mathcal{M}_B(\mO \otimes \Znl)
\end{equation}
\begin{equation}
F(N)=N^{\mathfrak{b}_q}
\end{equation}
Since taking $\mathfrak{b}_q$-invariants is exact on
$\mathfrak{B}_q$-modules and is faithful on $\MBG$ we get:
\begin{Prop}
This functor is exact and faithful.
\end{Prop}
The functor $F$ has a left adjoint. Let`s denote it $G$. By usual
Barr-Beck type theorem we get:
\begin{Prop}
 $F(G(\mO \otimes \Znl))=(\Dl)^{\mathfrak{b}_q}$ is a sheaf of algebras over $\Tl$.
 The category of modules over it is equivalent to $\DBG$.
\end{Prop}
We want to prove the following:
\begin{Prop}\label{equality}
$(\Dl)^{\mathfrak{b}_q}=D$ as sheaves of algebras.
\end{Prop}
\noindent \textit{Proof of proposition \ref{equality}.}
We have that 
\begin{equation}
(\Dl)^{\mathfrak{b}_q}=(O_q \otimes M_\lambda)^{\mathfrak{b}_q}=Ind_{\mathfrak{b}_q}^{G_q}(M_\lambda)=
Ind_{\mathfrak{u}_q}^{G_q}\circ Ind_{\mathfrak{b}_q}^{\mathfrak{u}_q}(M_\lambda)
\end{equation} 
hence it is enough to prove that
\begin{equation}  
  Ind_{\mathfrak{b}_q}^{\mathfrak{u}_q}(M_\lambda)=End_{\Znl} (M_{\lambda+2\rho}) 
\end{equation}  
  and that this map is comaptible with all relevant structures. Here we use two conventions about Verma modules, where $M_\lambda$ gets its action from the adjoint action of the restricted quantum group on the nonrestricted and in $End_{\Znl} (M_{\lambda+2\rho})$ we think of the Verma module as induced from the Borel. It is enough to prove that
\begin{equation}   
   Ind_{\mathfrak{b}_q}^{\mathfrak{u}_q}(M_{\lambda,\chi})=End (M_{\lambda+2\rho,\chi})
\end{equation}    

First note that if we denote by $o_q(G)$ the functions on the quantum Frobenius kernel (the fiber over the identity of sheaf of algebras $O_q(G)$ over $O(G)$) we have that $o_q(G)=o_q(N_-)\otimes o_q(B)$ and so  $Ind_{\mathfrak{b}_q}^{\mathfrak{u}_q}(M_{\lambda,\chi})=o_q(N_-)\otimes M_{\lambda,\chi}$ as vector spaces (actually as $N_{_q}$ modules). Now by definition $M_{\lambda,\chi}$ is an algebra (it is a quotient of $U_q(b)$ and using duality one can define an action of $o_q(N_-)$ on $M_{\lambda,\chi}$ as in \cite{J} getting the desired isomorphism by specializing the map from \cite{J} to the $\chi$ central character.


\subsubsection{Category $\DBG$ as a gluing of module categories
over quantum Weyl algebras.}
In this section we will give another way of viewing $\Dl$ as a sheaf of algebras over the cotangent bundle.

For $w \in \W$ De-Concini and Lyubashenko $\cite{DL}$ introduced
localizations $\mO_{q,w}$ of $\Oq$. These induce localizations
$\mO_w$ of $\mO$ corresponding to covering $G$ by translates of
the big cell $B_-B$. Joseph has introduced localizations of the
representation ring $\mO^{N_q}_q$ (\cite{J}), and it is easy to
see that they are exactly $\mO_{q,w}^{N_q}$. We thus have a
covering of the category of $\Oq$ modules by the categories
$\mO_{q,w}$. In other words
\begin{equation}\label{Gglue}
\Oq-mod = \underset{\longleftarrow}{\operatorname{lim}}_{w \in \W}
\mO_{q,w}-mod
\end{equation}

We have corresponding localizations $\MBG_w$ (i.e.
$B_q$-equivariant $\mO_{q,w}$-modules), of the category $\MBG$.
Then
\begin{equation}\label{MBGglue}
\MBG = \underset{\longleftarrow}{\operatorname{lim}}_{w \in \W}
\MBG_w
\end{equation}
Using the description of $\MBG$ as a Proj-category from \cite{BK},
it is clear that $\MBG_w$ is affine, i.e., $\mO_{q,w}$ is a
projective generator of $\MBG_w$ and hence $\MBG_w \cong
\mood-\End_{\MBG}(\mO_{q,w})$.The functors that induce this
equivalence are the adjoint pair
$(?\otimes_{\End_{\MBG}(\mO_{q,w})} \mO_{q,w}, ?^{B_q})$. Now
$\End_{\MBG}(\mO_{q,w}) = \Hom_{\MBG}(\Oq, \mO_{q,w}) =
\Gamma(\mO_{q,w})$. Joseph \cite{J} showed that:
\begin{equation}\label{JosephO}
\End_{\MBG}(\mO_{q,w}) = (M^\star_0)^{T_w}
\end{equation}
as $U_q$-modules, where $T_w$ are the algebra automorphisms of
$U_q$ from section 2.1.2. Note that these algebras are
non-isomorphic for different $w$ in general. When $w = e$, we have
$\End_{\MBG}(\mO_e) \cong U_q(\n_-)$ as an algebra.

Similarly to $\MBG$, $\MBG_w$ also has an induction functor from
$B_q$ modules, $V \mapsto \mO_{q,w} \otimes V$. It follows from
the fact that $\MBG_w$ is affine that
\begin{equation}
({\mO_{q,w} \otimes V})^{B_q} \cong \ {\mO_{q,w}^{B_q}} \otimes V
\end{equation}

 Next we have the forgetful functor $f_\star: \DBG \to \MBG$, which is exact and
faithful. It is easy to see that it has 
 left adjoint $f^\star: \MBG \to \DBG$ (just tensor with $\D$ over $\Oq$ and
factor out the necessary relations).

Define ${\DBG}_w$ to be the localization of $\DBG$ lying over
$\MBG_w$ (i.e., replace $\Oq$ by $\mO_{q,w}$ in the definition of
$\DBG$). We get adjoint pair of functors $(f^\star_w,
f_{w,\star})$ between these categories with the same properties as
above. Hence, abstract nonsense shows that
\begin{Prop}
$f^\star_w \mO_{q,w} = \mO_{q,w}\otimes_{\Oq} \Dl$ is a projective
generator of ${\DBG}_w$ and therefore \begin{equation} {\DBG}_w
\cong \mood-\End_{\DBG}(\mO_{q,w}\otimes_{\Oq} \Dl)
\end{equation}.
\end{Prop}

 Put
 \begin{equation}
 \Aw = \Gamma
(\mO_{q,w}\otimes_{\Oq} \Dl) =\End_{\DBG}(\mO_{q,w}\otimes_{\Oq}
\Dl)
\end{equation}.

 We can calculate explicitly these $\Aw$. We have
 \begin{equation}
 \Aw = \Gamma
(\mO_{q,w}\otimes_{\Oq} \Dl) = (\mO_{q,w}\otimes M_\lambda)^{B_q}
= \mO_{q,w}^{B_q}\otimes M_\lambda
\end{equation}
where the ring structure is induced from the one on $\D$ .

It follows that $\Aw$ coincide with the rings introduced by Joseph
\cite{J}. In the generic case these can be described as the
$B^{T_w}_q$-finite part of the full endomorphism ring
$\End((M^\star_{\lambda })^{T_w})$.

 When $w = e$  Joseph
describes explicitly this ring, a quantum Weyl algebra
\begin{equation}\label{qWalgebra}
\Ae \cong U_q(\n_-) \otimes U_q(\n_+)
\end{equation}
Here, the algebra structures on $U_q(\n_-)$ and on $U_q(\n_+)$ are
the usual ones and we have the commutation relations
\begin{equation}\label{commrelations}
q^{-(\alpha,\beta)}E_\alpha \otimes F_{\beta} - q^{(\alpha,
\beta)} F_\beta \otimes E_\alpha = \delta_{\alpha, \beta}
\end{equation}
for $\alpha, \beta \in \Delta$. (So $\Ae$'s algebra structure is
independent of $\lambda$.)

The category $\DBG$ can now be described as the gluing of module
categories
\begin{equation}\label{DBGglue}
\DBG = \underset{\longleftarrow}{\operatorname{lim}}_{w \in \W}
\Aw-\mood
\end{equation}

From the description of $\Aw$ as $\mO_{q,w}^{B_q}\otimes
M_\lambda$ with the algebra structure induced from $\D$ we get
\begin{Prop}\label{blablabla} $i)$ $Z(\Aw) = (\mO_w
\otimes \Znl)^B =\mO_w^B \otimes \Znl$ for each $w \in \W$. $ii)$
Consequently, $$qcoh(\Tl) = \underset{\longleftarrow}{\operatorname{lim}}_{w \in \W}
Z(\Aw)-\mood
$$ where $Z(\Aw)$ denotes the center of $\Aw$.
\end{Prop}

Let us remark that for generic $\lambda$ the variety $(G\times
N^{\lambda^{2l}})/B$ is an affine variety.

\subsubsection{Torsors} We first define a category $\tDBG$ of
$\D$-modules that contains all $\Dl$, $\lambda \in T_P$, a
"torsor".
\begin{defi}\label{d1t} An object of $\tDBG$
is a triple $(M, \alpha, \beta)$, where $\alpha: \D \otimes M \to
M$ a left $\D$-action and $\beta^{res}: M \to M\otimes \Oq(B)$ a
right $\Oq(B)$-coaction.
\smallskip

\noindent $i)$ The $U_q(\n_+)$-actions on $M$ given by $\beta\vert
U_q(\n_+)$ and by $(\alpha\vert_{U_q(\n_+)})$ coincide.

\noindent $ii)$ The map $\alpha$ is $U_q(\bb)$-linear with respect
to the $\beta$-action  on $M$ and  the action on $\D$.
\end{defi}

Put $\tM = U_q(\g)/\sum_{\alpha \in R_+} U_q(\g) \cdot E_\alpha$
(a "universal" Verma module) and define
\begin{equation}\label{tildeD}
\tD = \Oq \otimes \tM
\end{equation}
$\tD$ inherits an $U^{res}_q(\bb)$-module structure from $\D$, so
$\tD$ is an object in $\tDBG$.

$\tM$ has the $U^{res}_q(\bb)$-sub module $\Zbl = \Zl/(\Zl \cap
U_q \cdot U_q(\n_+)_+)$, which is an algebra generated by
$K^{(l)}_\mu, E^{(l)}_\mu,\, \mu \in R_+$. $b_q$ acts trivially on
$\Zbl$, so the $U^{res}_q(\bb)$-action factors to a
$U(\bb)$-action on $\Zbl$. We have:

\begin{Lem}\label{Bmod} As a
$B$-module $\Zbl \cong \mO(B)$ where the module structure comes as
follows: consider the map $B \to B$ defined as $tn \to t^2n$, this
is an unramified covering. Use this covering to pull back the
adjoint action of $B$ on itself .
\end{Lem}
\smallskip

\noindent \textit{Proof of proposition \ref{Bmod}.} To calculate
this one should notice that in the generic case the action induced
from the adjoint action of $U_q(n_+)$ on the universal verma
module is the same as the one induced from the commutator action.
Hence one can use the calculations of $\cite{CKP}$ which give the
required result. As for the torus part, the adjiont action is by
the grading and it is easy to see we get the required action.
$\Box$

 As in the previous sections we see
that $\mO(G \times B)$ imbeds to $\tD$, $B_q$-equivariantly. Note
that this embedding corresponds to the surjection ${(\;)}^{l}: T_P
\to T$. Define $B^{twist}$ to be the cover of $B$ induced from the
map $B= TN \ni tn \to t^{2l} n \in B$ with the natural
$B$-structure. We have
\begin{DT}\label{CT1y} We define $\tTX^{twist} = (G \times B^{twist})/B$ (where
$B$ acts on $G \times B^{twist}$ by $b \cdot (g,b') = (bg, b \cdot
b')$). $\tTX^{twist}$ is a $T_P$-torsor over $T^\star X$. $\tD$ is
a sheaf of Azumaya algebras over $\tTX^{twist}$.
\end{DT}

\subsection{Third construction and Azumaya splitting over fibers.}

Recall the $\UAres$ action on $U_q$. The center of $U_q$ is exactly the submodule of invariants with respect to the subalgebra $u_q$. Hence the center of $U_q$ is a $\mathfrak{g}$ module. Recall also that $U^{fin}$ is the finite part of $U_q$ with respect to the $\UAres$ action. We have:

\begin{Lem}\label{fincenter}
Let $\Z^{fin}$ be the center of $U^{fin}$. We have $\Z^{fin}=\Z \cap U^{fin}$. $Z^{fin}$ is an integrable  $\mathfrak{g}$ module and $Spec(\Z^{fin}$ is isomorphic to an unramified cover of $G \times _{T/W} T$ where the map $T \to T/W$ is induced from the l-th power map and the action is the adjoint action.
\end{Lem}
\noindent \textit{proof of proposition \ref{fincenter}} The first claim follows from Joseph`s \cite{J1} description of $\U$ inside $U_q$. The second statement follows from the calculations of \cite{DKP}.

Consider the diagram:

\begin{equation}\label{basic4}
\begin{matrix}
{T_P} & {\overset{}{\longleftarrow}} &
{\tTX^{twist}} \\
{{\downarrow}_{2l}} & {} & \downarrow_{{\pi}} \\
T/(\bullet, W) & \overset{}{\longleftarrow} & G
\end{matrix}
\end{equation}



\begin{defi}
$t \in T_P$ is unramified if for every $\alpha \in \tilde{\Delta}
= \Delta \cup \{-\alpha_0\}$, (where $\alpha_0$ is the longest
root), $(t(K_\alpha))^{2l} = 1$ implies $(t(K_\alpha))^2 =
q^{-(2\rho, \alpha)}$. $T^{unram}_P$ denotes the set of unramified
$t$'s.
\end{defi}
Note that any $t \in T_P$ can be made unramified by adding an
integral weight to it. We put $\Z^{unram} = \mO(T^{unram}_P
\times_{T/W} B_- B)$ and $U^{unram}_q = U_q \otimes_{\ZHC}
\mO(T^{unram}_P)$. Brown and Gordon \cite{BG} proved that
\begin{Lem}\label{BrownandGordon}
$i)$ $U^{unram}_q$ is Azumaya over $\Z^{unram}$ and $ii)$ for each
$\chi \in maxspec (\Z^{unram})$, we have $\End(M^{baby}_\chi) =
U_q/ \chi U_q$ (where $M^{baby}_\chi$ is the baby Verma module.)
\end{Lem}

 Let $\sigma: {\tTX^{twist}} \to T_P \times_{T/W} G$ be the
natural map. Note that this is a $G_q$-equivariant map. Now since both sheaves are also $G_q$-equivariant, from the description of the $\D$ as induced from
endomorphism of baby Vermas and from lemma \ref{BrownandGordon} that describes the sheaf over a dense subset (the big cell) we have

\begin{Prop}\label{Azsplit} The action map $U^{fin}_q \otimes_Z O_{\tTX^{twist}} \to
\tD $ induces an isomorphism $U^{fin}_q \otimes_Z
O_{\tTX^{twist,unram}} \cong \tD\mid _{\tTX^{twist,unram}} $.
\end{Prop}

Thus we see that over the preimage of the bigcell $\Dl$ is Azumaya. (Note again that in a formal neighbourhood of $p$ it is Azumaya everywhere).

Hence it follows as in \cite{BMR}(vanishing of the Brauer group of
a local ring with separably closed residue field) that $\tD$
Azumaya splits over the formal neighborhood of any fiber of
$\sigma$ lying over the big cell. Hence, we get

\begin{Cor}\label{AZTX} The category of $\Dl$-modules supported on
the fiber of $\sigma$ over the big cell is equivalent to $\mO$-modules supported on
the same fiber.
\end{Cor}

\subsection{Derived $\mathcal{D}$-affininty.}

\subsubsection{Global sections and vanishing cohomology of $\Dl$.}

\begin{Prop}\label{RGammaD} We have $i)$ $\Ul \cong R\Gamma(\Dl)$
and $ii)$ $\tU \cong R\Gamma(\tD)$ (if $l$ is a prime $p >$
Coxeter number of $G$.)
\end{Prop}
We recall from \cite{BK} that the natural map $\Ul \to
\Gamma(\Dl)$ was given as follows: There is the natural surjection
$\Ul \to M_\lambda$. It induces a surjective map
\begin{equation}\label{Inds1}
 \Ind \Ul \to \Ind M_\lambda = \Dl
\end{equation}
Since $\Ul$ is a $G_q$-module, \ref{Gmod} shows that $\Gamma(\Ind
\Ul) = \Ul$, which gives the desired map, by applying $\Gamma$ to
\ref{Inds1}.
\smallskip
\begin{Rem}\label{addgrade} Recall that $\Ul = U^{fin}/J_\lambda$.
In order to get global differential operators equal to
$U_q/J_\lambda$ one can enlarge $\Dl$ by adding some grading
operators; the vanishing of higher self extensions will remain
true.
\end{Rem}
 \noindent \textit{Proof of proposition \ref{RGammaD}.}
We have $\Dl = \Ind M_\lambda$. We have the integral version $\Ind
M_{\lambda, \A} \in \MBG_\A$.

Consider the specialization $\A \to \Fp$, $\nu \to 1$. In
\cite{BK} we showed that the statement about global sections in
$i)$ holds for a generic $q$. The argument given there transforms
to the case of $p$'th roots of unity if we can show that
\begin{equation}\label{vessla2}
\dim_\C \Gamma(\Ind \gr_i M_\lambda) = \dim_\C \Gamma (\Ind \gr_i
M_{\lambda, q=1})
\end{equation}
Now, by \cite{AJ}, $\dim_\C \Gamma (\Ind \gr_i M_{\lambda, q=1}) =
\dim_{\Fp} \Gamma (\Ind \gr_i M_{\lambda, \Fp})$. On the other
hand, it follows from \cite{APW}, that $\dim_\C \Gamma(\Ind \gr_i
M_\lambda) = \dim_{\Fp} \Gamma (\Ind \gr_i M_{\lambda, \Fp})$.

In order to prove that higher cohomologies vanishes in $i)$, it
suffices by \cite{APW} (page 26) to show that
\begin{equation}\label{vessla}
R\Gamma^{>0}_{\Fp}(\Ind M_{\lambda, \Fp}) = 0
\end{equation}
This holds, because $R\Gamma^{>0}(\Ind \gr M_{\lambda, \Fp}) =
R\Gamma^{>0}(\Ind S(\n_{-,p})) = 0$, by \cite{AJ}. ($\gr$ is taken
with respect to the filtration on $M_\lambda$ coming from the
identification $M_\lambda = U_q(\n_-)$ and putting each $F_\mu$ in
degree $1$.) This proves $i)$. $ii)$ is similar. $\Box$

Now as in \cite{BK} we can extend this to almost all roots of unity.

\subsubsection{Localization functor.} Recall the definition of the
localization functor $Loc_\lambda: \Ul-\mood \to \DBG$ \cite{BK}. 
\begin{defi}\label{localization}
Define the localization functor
$$Loc_\lambda: \Gamma(\Dl)-\operatorname{mod} \to \DBG$$ by
$M \to \Dl \otimes_{\Ul} M$, where we have used $\Ul = \Gamma(\Dl)$ .
\end{defi}
This is a left adjoint to the global sections functor. Note that
$Loc_\lambda(\Ul)=\Dl$.
 Similarly, we can define localization $\tLoc: U_q^{fin}-\mood \to
\tDBG$. $Loc_\lambda$ has a left derived functor $\Ll:
D^b(\Ul-\mood) \to D^b(\DBG)$ and $\tLoc$ has left derived functor
$\tL:D^b(U_q^{fin}-\mood) \to \tDBG$.

\begin{Lem}\label{adjoint} $\tL$ is left
adjoint to $R\widetilde{\Gamma}$ and $\Ll$ is left adjoint to
$R\Gamma$ (if $\lambda$) is regular.
\end{Lem}
If $\lambda$ isn't regular the second statement in the lemma will
hold if we replace $D^b$ by $D^{-}$.
\begin{Prop}\label{closetoequi} $i)$ The functor $R\widetilde{\Gamma} \circ
\tL: D^b(U_q) \to D^b(\tU)$ is isomorphic to the functor $M \to M
\otimes_{\ZHC} \mO(T_P)$. $ii)$ For regular $\lambda$, the
adjunction map $id \to R\Gamma \circ \Ll$ is an isomorphism.
\end{Prop}
\noindent \textit{Proof.} $i)$  follows from part $ii)$ of lemma
\ref{RGammaD} for free modules and then from the same lemma again
for general modules, by considering free resolutions.

For $ii)$ observe that for regular $\lambda$, for any $M \in
D^b(\Ul-\mood)$, we have canonically $M \otimes_{\ZHC} \mO(T_P) =
\oplus_{w \in \mathcal{W}}\oplus_{\Hom(P, \{+1, -1\})} M$. Now the
claim follows since $R\Gamma \circ \Ll(M)$ is one of these direct
summands. $\Box$ 

\subsubsection{Derived $\mathcal{D}$-affinity.}
\begin{Thm}\label{Daffthm}
$R\Gamma: D^b(\DBG) \to D^b(\Ul-\mood)$ is an equivalence of
categories.
\end{Thm}
Noting that the canonical bundle of $\Tl$ is trivial and that
$\pi$ from \ref{basic4} is a projective morphism the theorem now
follows from the following lemma which is a slight generalization of \cite{BMR}:
\begin{Lem}\label{keykey}
Let $\A$ be a generically Azumaya algebra over a smooth variety $X$(i.e. Azumaya over a generic point). Suppose
that $X$ is Calabi-Yau (i.e. $\omega_X \cong \mO_X$) and that we
have a projective map $\pi: X \to Spec R$ for some commutative
algebra $R$. Suppose also that the derived global section functor
$R\Gamma: D^b(\A-\mood) \to D^b(\Gamma(\A)-\mood)$ has a right
adjoint $\mathcal{L}$ and the adjunction morphism $id \to R\Gamma
\circ \mathcal{L}$ is an isomorphism. Then $R\Gamma$ is an
equivalence of categories.
\end{Lem}

Note that, for generic $\lambda$, $\Tl$ is affine and hence we
then get an equivalence $\Gamma:\DBG \cong \Ul-\mood$.

\section{Applications} Assume for simplicity $\lambda$ is integral
and regular and that $\chi$ belongs to $B_- B$ and is unipotent.
We know that
\begin{equation}\label{final1}
D^b(\DBG) \cong D^b_\lambda(U_q-\mood)
\end{equation}
We get that
\begin{equation}\label{final2}
D^b(\DBG_{(\chi, \lambda)}) \cong D^b_{\hat{\chi},\lambda}(U_q-\mood)
\end{equation}
where the left hand side denote those (complexes of) $\Dl$-modules
supported on the Springer fiber of $(\chi,\lambda)$ and the right
hand side denote $\Ul$ modules who are locally annihilated by a
power of the maximal ideal in $\Z$ corresponding to $\chi$
(generalized central $l$-character $\chi$).

By the Azumaya splitting we have
\begin{equation}\label{final3}
\DBG_{(\chi, \lambda)} \cong qcoh(T^\star X)_{\chi}
\end{equation}
where the latter category is the quasi coherent $\mO$-modules on
$T^\star X$ supported on the Springer fiber of $\chi$ (with
respect to the usual Springer resolution $T^\star X \to
\mathcal{N} =$ unipotent variety of $G$). Note that for the trivial central character this equivalence is Koszul dual to the equivalence in \cite{ABG} thus giving a geometric proof of their equivalence.

We deduce an equivalence of $K$-groups
\begin{equation}\label{final4}
K_{\hat{\chi}}(\Ul-\mood) \cong K(qcoh(T^\star X)_{\chi})
\end{equation}
Note that $K(qcoh(T^\star X)_{(\chi)})$ is torsion-free. From
\ref{final4} we see for instance that the number of irreducible
$U_{q,\chi}$-modules equals the rank of $K(qcoh(T^\star
X)_{(\chi)})$.

In the following paper we will use the methods here to prove a conjecture made by DeConcini, Kac and Processi stating that the dimensions of irreducible $U_{q,\chi}$-modules is divisible by $l^{dim(O)/2}$
where $l$ is the order of the root of unity and $O$ is the orbit through $\chi$.  

Notice also that \cite{BMR} showed that $K(qcoh(T^\star X)_{\chi})
\cong K(qcoh(T^\star X)_{\chi}(\bar{\mathbb{F}}_p))$ for $p >$
Coxeter number of $G$. This relates our work to the representation
theory of $\g_p$ via the results in \cite{BMR}.

Actually, it is possible to relate our work with \cite{BMR}
as to get more transparent geometric proofs of the theorems
appearing in the work of Andersen, Jantzen and Soergel \cite{AJS},
leading to the proof of Lusztig`s conjecture about multiplicities
in characteristic p and also extending them to nontrivial central characters. Our work defines a perverse t-structure on O
modules in zero characteristic (root of unity), where \cite{BMR}
define such a structure in positive characteristic. By showing
that specializing the root of unity t-structure one would get (at
least for big p) the positive characteristic t-structure, one
would be able to deduce that the quantum modules and the positive
characteristic modules have the same multiplicity formulas, as
shown in \cite{AJS}. This will appear in future work.


\begin{thebibliography}{}
\bibitem[AG]{AG} S. Arkhipov, D. Gaitsgory,
\textit {Another realization of the category of modules over the small quantum group.}
  Adv. Math.  173  (2003),  no. 1, 114--143.

\bibitem[AJ]{AJ} H. H. Andersen, J. Jantzen,
\textit{Cohomology of induced representations for algebraic
groups}, Math. Ann. 269, (1984) 487-525.

\bibitem[AJS]{AJS} H. H. Andersen, J. Jantzen, W. Soergel,
\textit{Representations of quantum groups at a pth root of unity
and of semisimple groups in characteristic p: independence of p},
Asterisque 220 (1994).

\bibitem[APW]{APW} H. H. Andersen, P. Polo and Wen Kexin,
\textit{Representations of quantum algebras}, Invent. Math. 104,
(1991) 1-59.




\bibitem[BK]{BK} E. Backelin and K. Kremnitzer
\textit{Quantum flag varieties, equivariant quantum
$\mathcal{D}$-modules, and localization of Quantum groups.}, to
appear in Adv. in Math., arXiv:math.QA/0401108.

\bibitem[BB]{BB} A. Beilinson and J. Bernstein,
\textit{Localisation de $\g$-modules}, C. R. Acad. Sc. Paris, 292
(S\'erie I) (1981) 15-18. 




\bibitem[BMR]{BMR} A. Bezrukavnikov, I. Mircovic and D. Rumynin,
\textit{Localization for a semi-simple Lie algebra in prime
characteristic.}, arXiv:math.RT/0205144.

\bibitem[BG]{BG} K. A. Brown, I. Gordon,
\textit{The ramification of centeres: Lie algebras in positive
characteristic and quantizied enveloping algebras.},
arXiv:math.RT/9911234

\bibitem[CP]{CP} W. Chari and A. Pressley,
\textit{A guide to quantum groups}, Cambridge University Press,
Cambridge 53, (1995).

\bibitem[DL]{DL} C. De Concini and V. Lyubashenko,
\textit{Quantum function algebras at roots of $1$}, Adv. in Math.
108, (1994) 205-262.


\bibitem[CKP]{CKP} C. De Concini, V.G. Kac and C. Procesi,
\textit{Quantum coadjoint action }, JAMS v.5, number 1, (1992)
151-189.





\bibitem[J]{J} A. Joseph,
\textit{Faithfully flat embeddings for minimal primitive quotients
of quantized enveloping algebras}, In A. joseph and S. Shnider
(eds.), Quantum deformation of algebras and their representations,
Israel Math. Conf. Proc. 7 (1993), pp 79-106.

\bibitem[JL]{JL} A. Joseph, G. Letzter
\textit{Local finiteness for the adjoint action for quantized
enveloping algebras}, J. Algebra 153, (1992), 289-318.




\bibitem[M]{M} S. Montgomery,
\textit{Hopf Algebras and Their Actions on Rings}, CBMS (1993).












\end{thebibliography}
\end{document}